\documentclass[12pt]{article}
\usepackage{epsfig}
\usepackage{amsmath}

\newcommand{\no}{\nonumber}

\newcommand{\be}{\begin{equation}}

\newcommand{\ee}{\end{equation}}

\newcommand{\ba}{\begin{eqnarray}}

\newcommand{\ea}{\end{eqnarray}}

\newcommand{\fee}{\mbox{$\varphi$}}

\newcommand{\lam}{\mbox{$\lambda$}}

\begin{document}
\title{\hspace{.2 in} \\ 
{\Large \bf A Many to One Discrete Auditory Transform}}

\author{Jack Xin \thanks{Corresponding author, Department of Mathematics,
UC Irvine, Irvine, CA 92697, USA; email:jxin@math.uci.edu.} \hspace{.1 in} and 
\hspace{.1 in} Yingyong Qi \thanks{
Qualcomm Inc,
5775 Morehouse Drive,
San Diego, CA 92121, USA.}
}

\date{}

\maketitle

\vspace{0.25 in}
\thispagestyle{empty}

\begin{abstract}
A many to one discrete auditory transform is presented to map a sound signal
to a perceptually meaningful spectrum on the scale of human auditory filter band widths 
(critical bands). A generalized inverse is constructed in 
closed analytical form, preserving the band energy and 
band signal to noise ratio of the input sound signal. The forward and inverse 
transforms can be implemented in real time. 
Experiments on speech and music segments 
show that the inversion gives a perceptually equivalent 
though mathematically different sound from the input. 

\end{abstract}

\vspace{.4 in}

\hspace{.2 in} {\bf Keywords: Many to One Discrete
Auditory Transform}
\vspace{ .1 in} 

            
\newpage
\setcounter{page}{1}
\section{Introduction}
\setcounter{equation}{0}
Short term discrete Fourier transform (DFT) is a common tool to 
map sound signals from time 
domain to spectral domain for analysis and synthesis \cite{Br}. However, the spectral resolution 
of DFT over a standard short time window of 5 to 15 milleseconds (ms) 
in duration is much more refined than 
the resolution of human auditory filters that have band widths referred to as critical bands \cite{Hart,ZF}.  
Critical bands are nearly uniform in widths similar to DFT for frequencies 
under 500 Hz, yet the widths increase rapidly towards higher frequencies. The 
nonuniform frequency resolution of the ear resembles that of wavelets \cite{ID,Strang}, 
though critical band widths do not follow a simple power law, and auditory 
filter shapes may not obey the requirements of the wavelet basis functions. 
An orthogonal discrete 
transform with broader and smoother spectrum towards higher frequencies than that of DFT
is recently constructed \cite{xq05} to mimic the auditory filtering. Due to the 
limitation of orthogonality, the variation of the spectrum does not match the scale  
of critical bands. In addition, the spectrum of the transform does not 
carry enough perceptual meaning and so makes it inconvenient to perform psychoacoustically 
based spectral analysis and processing. 
\medskip

In this paper, we present a novel many-to-one discrete auditory transform (MDAT) that 
maps sound signals from the time domain to a perceptually meaningful spectral domain
on the scale of critical bands. 
The many-to-one mapping is consistent with the fact that physically and 
mathematically different signals can sound the same to human ears \cite{Hart,Sch,ZF}. 
The frequency resolution for the perception of sound in our brain is much 
lower than that is required to fully describe a signal mathematically \cite{Sch}. 
The perception variables of MDAT are band energies and band signal to noise ratios (SNRs), 
motivated by perceptual coding in AAC and MP3 technology of digital music compression 
\cite{AAC,Poh}. The SNRs depend on two neighboring frames of a signal and 
so MDAT spectrum also encodes temporal information, different from DFT. 
As a test of the efficiency of these variables, and for the synthesis of sounds 
post spectral processing, we show how to construct an inverse which is perceptually 
equivalent to the input sound though mathematically not identical. Both the forward and 
inverse operations are in closed analytical form, and allow real time implementation 
of the resulting algorithms.   
\medskip

Compared with DFT (implemented by FFT), MDAT has better temporal resolution 
due to its lower spectral resolution in higher frequencies. 
In terms of the 256 point FFT used in this 
paper, the number of frequency bands of MDAT 
is in the 40's (see Tables 1 and 2 for signals with different sampling frequencies) while 
DFT has 128 frequency components. 
Compared to time domain filter bank with a
relatively small number of band pass filters (4 to 16 channels) as in body-worn hearing devices \cite{Green}, 
MDAT has better frequency resolution yet does not have the delays encountered when 
a larger number of frequency separating band pass filters are needed. Hence MDAT is 
expected to be a useful tool in 
applications where a spectral processing strategy is necessary on the critical band 
scale, and a trade-off of spectral accuracy and temporal precision is to be optimized.

\medskip

The paper is organized as follows. In section 2, the MDAT is formulated and the 
associated perceptual variables are defined. 
Then an inverse is constructed in closed analytical form based on band energies 
and SNRs. In section 3, MDAT is applied to speech (sampled at 16 kHz) and music (sampled at 44.1 kHz) signals, and 
properties of perceptual spectral variables are illustrated. 
The reconstructed signals are compared with the input signals both spectrally 
and in waveforms, and these signals can be heard at author's website \cite{web}. 
Section 4 contains discussion and conclusion. 

\section{MDAT and a Perceptual Inversion}
\setcounter{equation}{0}
Let $s=(s_0,\cdots,s_{N-1})$ be a discrete real signal, 
the discrete Fourier transform (DFT) is \cite{Br}:
\be
 \hat{s}_{k} = \sum_{n=0}^{N-1}\, s_{n}\, e^{-i (2\pi n k/N )}. \label{p0}
\ee
The DFT is implemented by the fast Fourier transform (FFT) algorithm, 
we shall refer to the $k=0$ component of DFT as DC (direct current) and 
the other components as AC (alternating current) for short.

Let us further map the $\hat{s}_{k}$'s to a spectral domain of lower resolution  
where perception variables can be better defined. Such a spectral domain is obtained from 
binning the DFT components into bands of various widths, similar to 
the critical band width distribution of human auditory filters. 
The detailed partition of DFT components, the band widths, and 
psychoacoustic bark values of the bands are listed in Table 1 and 
Table 2. Table 1 is at sampling frequency $Fs=16$ kHz for speech sounds, 
and Table 2 is at $Fs=44.1$ kHz for music sounds.  
 Let $b = 1,2,\cdots, J$, denote the number of bands,  
and let $B(b)$ denote the 
DFT wave numbers $k$ in the $b$-th band. In case of Table 1, $N/J \approx 5.56$; and 
in Table 2, $N/J \approx 6.24$. 

The signal energy in the $b$-th band is:
\be 
e(b) = \sum_{k \in B(b)}\, |\hat{s}_{k}|^{2}. \label{p1a} 
\ee
Let $snr^b$ be the signal to noise ratio (SNR) in the b-th band, the perception domain
consists of nonnegative $2J$-dimensional vectors whose components are band energies and band SNRs:
\be
V_{perc}=\{ ( e(1),snr^{1}, e(2),snr^{2},\cdots, e(J), snr^{J})\}. \label{p1}
\ee
The $snr^b$ are calculated following the AAC coding \cite{AAC}, 
an improvement of MP3 coding \cite{Poh}. Let $r(k,t)$ and $f(k,t)$ be the amplitude and phase of 
$\hat{s}(k)$ at time frame $t$ denoted by $\hat{s}(k,t)$. 
The predicted amplitude and phase at time frame $t$ are:
\ba
r_{pred}(k,t) = r(k,t-1) + \Delta r, \; \Delta r \equiv r(k,t-1)-r(k,t-2), \no\\
f_{pred}(k,t) = f(k,t-1) + \Delta f,\; \Delta f \equiv f(k,t-1)-f(k,t-2). \label{p2}
\ea
The unpredictability measure of the signal, a quantity for 
measuring the noisy (uncertain) part of signal, is:
\be
c(k,t)={{\rm abs}(\hat{s}(k,t) - \hat{s}_{pred}(k,t)) \over {\rm abs}(\hat{s}(k,t))+{\rm abs}(\hat{s}_{pred}(k,t))},
\label{p3}
\ee
where $\hat{s}_{pred}(k,t) = r_{pred}(t)\, e^{i f_{pred}(t)}$. 
It is clear that $c(k,t) \in [0,1]$. Note that $c(k,t)$ encodes the 
time domain information of the signal $s$, which is not available in DFT. 
As a result, the perceptual variables (\ref{p1}) has both spectral and temporal 
information of the input signal. 
We shall omit the $t$ dependence from now on, 
as all subsequent operations will not explicitly use $t$.

The weighted unpredictability measure is:
\be
ec(b) = \sum_{k \in B(b)}\, r^{2}(k)\, c(k). \label{p4}
\ee
Next, convolve $e(b)$ and $ec(b)$ with spreading functions \cite{Sch} 
on the bark scale \cite{Hart} as:
\ba
ecb (b)& = & \sum_{b'=1}^{J}\, e(b')\, {\rm spread}\, ({\rm bark}\, (b'),{\rm bark}\,(b)), \label{p5} \\
ct(b) & = & \sum_{b'=1}^{J} \, ec(b')\, {\rm spread}\, ({\rm bark}\, (b'),{\rm bark}\,(b)), \label{p6}
\ea
where bark$(b)$ is the bark value of the b-th partition (band). The bark scale \cite{Hart} 
is nearly uniform on the logarithmic frequency scale. The spreading functions 
\cite{Sch} carry the shape information of human auditory filters.

Normalizing $ct$ by energy $ecb$ gives: 
\be
cb(b) = ct(b)/ecb(b), \label{p7}
\ee
a noise to signal ratio, which in turn defines tonality index as:
\be
tb(b)= -0.299 - 0.43\, \log (cb(b)), \label{p8}
\ee
if the value is in $(0,1)$, otherwise equal to zero if the value is below zero, or one if the value is 
above 1. Finally, the signal to noise ratio in decibel (dB) is:
\be
snr^b = tb(b)\, {\rm TMN} + (1 - tb(b))\, {\rm NMT}, \label{p9}
\ee
where ${\rm TMN} = 18$ dB (tone masking noise), ${\rm NMT}=6$ dB (noise masking tone). 
The forward transform denoted by $T$ from signal $s$ to its image in the perception 
domain $V_{perc}$ is a many-to-one mapping. Clearly, $T (-s) = T s$. 
\bigskip

We notice that each $snr^{b}$ is a monotone function of $cb(b)$ which in turn depends on 
$ec(b)$ and $ct(b)$. So two other ways of characterizing the perception domain are:
\be
V_{perc}^{(1)}=\{ ( e(1), cb(1), e(2), cb(2),\cdots, e(J), cb(J))\}, \label{p10a}
\ee
\be
V_{perc}^{(2)}=\{ (e(1), ec(1), e(2), ec(2),\cdots, e(J), ec(J))\}. \label{p10b}
\ee
In other words, $V_{perc}^{(1)}$ or $V_{perc}^{(2)}$ is sufficient to describe 
the perception variables, i.e. the band energies and band SNRs. 
Below we show how to reconstruct a sound signal 
from $V_{perc}^{(1)}$ or $V_{perc}^{(2)}$ and obtain a perceptually equivalent 
inverse.

\bigskip

The inversion from a subset of $2J$ dimensional space 
to the signal space $R^N$ ($N > 2J$) is non-unique. 
The inversion is through reconstructing the DFT vector $\hat{s}_{k}$. 
Let us write the reconstructed DFT vector as:
\be
a^{b}_{k} = w^{b}_{k}\, e^{1/2}(b)\, e^{i \fee^{b}_{k}}, \;\; k \in B(b), \label{p10}
\ee
where the real weighting factors $w^{b}_{k}$ satisfy for all $b$:
\be
 \sum_{k \in B(b)}\, |w^{b}_{k}|^2 = 1, \label{p11}
\ee
to preserve the band energy $e(b)$. The real phase factors $\fee_{k}^{b}$, and the DC component of 
DFT are assumed to be known for the reconstruction of the AC part of the DFT amplitude.  

\bigskip

The second conserved quantity (constraint) is $ec(b)$ in (\ref{p4}):
\be 
ec(b)  =  \sum_{k \in B(b)}\, |w^{b}_{k}|^2 \, e(b)\, c(k) 
=  e(b)\, \sum_{k \in B(b)}\, |w^{b}_{k}|^2 \, c(k). \label{p11a}
\ee
Define:
\be
<w^b>_{c}^2 = \sum_{k\in B(b)}\, |w^{b}_{k}|^2 \, c(k), \label{p11b}
\ee
which equals
\be
<w^b>_{c}^2 = ec(b)/e(b) \in (\min_{k \in B(b)}\, c(k), \max_{k \in B(b)}\, c(k)) \subset [0,1]. 
\label{p12}
\ee

If the inversion is from $V^{(2)}_{perc}$, then the two spectral constraints (\ref{p11}) 
and (\ref{p11b}) are available to be imposed in each band 
containing at least two DFT components. If the inversion is from $V^{(1)}_{perp}$, then 
$<w^b>_{c}^{2}$ has to be recovered from $e(b)$ and $cb(b)$.  
By (\ref{p7}), we have for each $b \in [1,J]$:
\be
cb(b) = {\sum_{b'=1}^{J}\, e(b')\, <w^{b'}>_{c}^{2}\, {\rm spread}({\rm bark}(b'),{\rm bark}(b))
\over \sum_{b'=1}^{J}\,e(b')\, {\rm spread}({\rm bark}(b'),{\rm bark}(b))}, \label{p13}
\ee
or
\ba
& & \sum_{b'=1}^{J}\, e(b')\, <w^{b'}>_{c}^{2}\, {\rm spread}({\rm bark}(b'),{\rm bark}(b)) \no\\
& = & 
cb(b) \sum_{b'=1}^{J}\, e(b')\, {\rm spread}({\rm bark}(b'),{\rm bark}(b)). \label{p14}
\ea
Equation (\ref{p14}) can be recast as a matrix equation $S \vec{x} = \vec{z}$, 
where $S = ({\rm spread}({\rm bark}(b'),{\rm bark}(b))$ 
is a square matrix, $\vec{x}$ is the column vector with entries $e(b)\, <w^b>_{c}^{2}$, 
$\vec{z}$ the right hand side column vector. 
The commonly used spreading matrix $S$ (based on e.g. Schroeder's spreading functions \cite{Sch}) 
does not have a nonnegative inverse. In order to find 
nonnegative solutions in general, one may solve a quadratic programming problem 
from (\ref{p13}). Define the matrix $Q=(q_{ij})$ with its entries:
\[ q_{ij}= { e(j)\, {\rm spread}({\rm bark}(j),{\rm bark}(i))
\over 
\sum_{j=1}^{J}\, e(j)\, {\rm spread}({\rm bark}(j),{\rm bark}(i))}. \]
The matrix $Q$ is invertible. A column vector $\vec{y}=(<w^b>_{c}^{2})$ is 
sought to minimize the $l^2$ norm $\| \vec{cb} - Q \vec{y}\|_{2}$ subject to the constraint 
$yl(b) \leq y(b) \leq yu(b)$, $yl(b)= \min_{k \in B(b)}\, c(k)$, $yu(b)=\max_{k \in B(b)}\, c(k)$. 

\medskip

In signal processing tasks that keep the band SNRs invariant as in hearing aids 
gain prescriptions, the quadratic programming is not needed, directly inverting 
$S$ will suffice to find $(<w^b>_{c}^{2})$. 

\medskip

Next we solve for $w^{b}_{k}$ from the two 
equations (\ref{p11}) and (\ref{p11b}), using information of $c(k)$, $k \in B(b)$. 
Let $N_b$ be the number of DFT components in $B(b)$, 
$\vec{\rho}= (|w^{b}_{k_1}|^{2}, |w^{b}_{k_2}|^{2},\cdots, |w^{b}_{k_{N_b}}|^{2})^{T}$, 
$\vec{\psi} = (c(k_1),c(k_2),\cdots, c(k_{N_b}))^{T}$, $k_j \in B(b)$, 
$\theta_{b} = <w^{b}>_{c}^{2}$, $\vec{e}=(1,1,\cdots,1)^{T} \in R^{N_b}$, $T$ denoting 
transpose. Equations (\ref{p11}) and (\ref{p11b}) now read (dot refers to inner product):
\ba
& & \vec{e}\cdot \, \vec{\rho}\, =\, 1, \label{p14a}\\
& & \vec{\psi}\cdot \, \vec{\rho} \, = \, \theta_{b}. \label{p14b}
\ea

If $\vec{\psi}$ is parallel to $\vec{e}$, equation (\ref{p14b}) is 
redundant with $\theta_b = c(1)$ by definition and equation (\ref{p14a}). 
This is true in particular if $N_b =1$.
The simplest smooth solution to (\ref{p14a}) is $\vec{\rho} = {1\over N_b} \vec{e}$.

If $N_b \geq 2$ and $\vec{\psi}$ is not parallel to $\vec{e}$, define vector:
\be
\vec{v} = \vec{e} - {\vec{e}\cdot \vec{e} \over \vec{e}\cdot \vec{\psi}}\, \vec{\psi}  \not = 0,
\label{p15}
\ee
clearly $\vec{v}\cdot \vec{e} = 0$, and $\vec{e}\cdot \vec{\psi} > 0$.
Equations (\ref{p14a}) and (\ref{p14b}) imply that:
\be
\vec{v}\cdot \vec{\rho} = 1 - 
{\vec{e}\cdot \vec{e} \over \vec{e}\cdot \vec{\psi}}\,\theta_b. \label{p16}
\ee

Equation (\ref{p14a}) and equation (\ref{p16}) say that in the orthonormal basis 
with $\vec{e}$ and $\vec{v}$ as two directions, the coordinates along $\vec{e}$ and $\vec{v}$ 
are constrained, the other coordinates are free. 
The simplest two dimensional solution is obtained 
by setting the free coordinates to zero ($\|\cdot \|_2$, $l^2$ norm or the Euclidean distance):
\be
\vec{\rho} = {1\over \sqrt{N_b}}\, {\vec{e}\over \sqrt{N_b}} + {1\over \|\vec{v}\|_2} 
(1 - {\vec{e}\cdot \vec{e} \over \vec{e}\cdot \vec{\psi}}\theta_b)\, {\vec{v} \over \| \vec{v}\|_2},
\label{p17}
\ee
which becomes upon substituting in (\ref{p15}):
\be
\vec{\rho} = \left [ {1\over N_b} + {1\over \|\vec{v}\|^{2}_{2}} 
\left (1 - {\vec{e}\cdot \vec{e} \over \vec{e}\cdot \vec{\psi}}\theta_b \right) \right ]\vec{e}
- {\vec{e}\cdot \vec{e} \over \vec{e}\cdot \vec{\psi}}\,{1\over \|\vec{v}\|^{2}_{2}} 
  \left (1 - {\vec{e}\cdot \vec{e} \over \vec{e}\cdot \vec{\psi}}\theta_b \right) \vec{\psi}.
\label{p18}
\ee
The regularity of solution (\ref{p17}) or (\ref{p18}) is no worse than that of 
$\vec{\psi}$ which is oscillatory in general. 
With the $w^{b}_{k}$'s so determined, a time domain signal is  
reconstructed by inverse DFT using the reconstructed 
$a^{b}_{k}$, $k \in B(b)$, $b =1, 2, \cdots, J$.

If $N_b=2$, (\ref{p18}) is the unique solution. If $N_b \geq 3$ (true if $b$ is above 
some critical number, see Table 1 and Table 2), there are infinitely many 
solutions to (\ref{p14a})-(\ref{p14b}). It is desirable to seek
a smoother solution because spectral smoothness improves temporal localization 
of the inverse transform. One way to obtain a smoother solution over the frequency bands 
($N_b\geq 3$) starting with FFT wave number $k_0$ is to minimize the following quadratic function:
\be
f = {1\over 2} \, \sum_{k=k_0}^{k_0+M-1} (\rho_{k+1}-\rho_{k})^{2}, \label{p19}
\ee
where $M+1$ is the total number of DFT components in those bands $B(b)$ with $N_b \geq 3$, 
subject to the two constraints (\ref{p14a})-(\ref{p14b}) in each such band $B(b)$.
Let $\vec{u}=(\rho_{k_0},\cdots, \rho_{k_0 + M})^{T}$, and define:
\ba
g_b(\vec{u}) & = & -1 + \sum_{k \in B(b)}\, \rho_{k}, \label{p20}\\
h_b(\vec{u}) & = & -\theta_b +\sum_{k \in B(b)}\, c_k \, \rho_k, \label{p21}
\ea
then the constraints are of the form $g_b= 0$ and $h_b = 0$. 
The minimizer can be approached as a steady state in  
a constrained gradient descent method \cite{Osher}. 
Let $\vec{u} =\vec{u}(\tau)$ solve the equation:
\be
\vec{u}_\tau = -\nabla_{\vec{u}}\, f 
-\sum_{b, N_b \geq 3}\, \lam_b \, \nabla_{\vec{u}}\, g_b - \sum_{b, N_b \geq 3}\, 
\eta_b \, \nabla_{\vec{u}}\, h_b, \label{p23}
\ee
where the Lagrange multipliers $\lam_b$ and $\eta_b$ are chosen so that the 
constraints in each band are satisfied:
\ba
{d\over d\tau}\, g_b(\vec{u}) & = & \nabla_{\vec{u}}\, g_b \cdot \vec{u}_{t} \no \\
& = & 
-\nabla_{\vec{u}}\, g_b \cdot \nabla_{\vec{u}}\, f  - 
\lam_b \, |\nabla_{\vec{u}}\, g_b|^{2} - 
\eta_b \nabla_{\vec{u}} g_{b}\cdot \nabla_{\vec{u}} h_b = 0, \label{p24}
\ea
\ba
{d\over d\tau}\, h_b(\vec{u}) & = & \nabla_{\vec{u}}\, h_b \cdot \vec{u}_{t} \no \\
& = & 
-\nabla_{\vec{u}}\, h_b \cdot \nabla_{\vec{u}}\, f  - 
\lam_b \, \nabla_{\vec{u}}\, h_b \cdot \nabla_{\vec{u}} g_b - 
\eta_b |\nabla_{\vec{u}} h_b |^{2} = 0. \label{p25}
\ea
We have used the fact that $\nabla_{\vec{u}} h_b$ or $\nabla_{\vec{u}} g_b$
only have nonzero components in the band $B(b)$. 
To solve (\ref{p24})-(\ref{p25}) band by band, it is convenient to consider 
\be
\tilde{c}_{k} = 1 - {c_k N_b \over \sum_{j \in B(b)}\, c_j }, \; k \in B(b). \label{p26a}
\ee
If $\tilde{c}_{k} = 0$, for all $k \in B(b)$, 
then the second constraint $h_b =0$ is redundant, $\eta_b=0$, and 
\be
\lam_b = -{ \nabla_{\vec{u}} g_b \cdot \nabla_{\vec{u}} f \over |\nabla_{\vec{u}} g_b|^2}.
\label{p26}
\ee
 
If $\tilde{c}_{k} \not = 0$, for some $k \in B(b)$, replace the constraint $h_b = 0$ by:
\be
\tilde{h}_{b}(\vec{u}) = \sum_{k \in B(b)}\, \tilde{c}_{k}\rho_k - 1 +  
{N_b <w^b>_{c}^{2} \over \sum_{k \in B(b)} \, c_k } = 0. \label{p27}
\ee
Then the ${\vec{u}}$ equation is (\ref{p23}) with $\tilde{h}_b$ in place of $h_b$. 
Due to $\nabla_{\vec{u}}\, g_b \cdot \nabla_{\vec{u}}\, \tilde{h}_{b} = 0$, 
$\lam_b $ is as given in (\ref{p26}), and:
\be
\eta_b = - { \nabla_{\vec{u}} \tilde{h}_b \cdot 
\nabla_{\vec{u}} f \over |\nabla_{\vec{u}} \tilde{h}_b|^2},
\label{p28}  
\ee
where $|\nabla_{\vec{u}} \tilde{h}_b|^2 = \sum_{k \in B(b)}\, \tilde{c}_{k}^{2}$, and
$|\nabla_{\vec{u}} g_b|^2 = N_b$ in (\ref{p26}). 

\bigskip

Finally, let us put the ${\vec{u}}$ equation in matrix form. 
Let $A$ be the symmetric tridiagonal matrix with $1$'s on the off-diagonals, and 
$(-1,-2,\cdots,-2,-1)$ on the diagonal ($\cdots$ refer to $-2$'s), 
then $\nabla_{\vec{u}}\, f = A \vec{u}$. 
Let $R$ be the block diagonal matrix where each block is the symmetric $N_b\times N_b$ 
matrix with the $(i,j)$-th entry being 
$N_{b}^{-1} + {\tilde{c}_{i}\tilde{c}_{j} \over \sum_{k \in B(b)}\, \tilde{c}_{k}^{2}}$.
If $\sum_{k \in B(b)}\, \tilde{c}_{k}^{2}$ is zero, 
the second term in the sum is understood to be absent. 
The matrix form of $\vec{u}$ equation is ($I$ the identity matrix):
$\vec{u}_{\tau} = (I - R)A \vec{u}$, 
whose solution is in closed form $\vec{u}(\tau) =\exp\{ (I-R)A \tau \} \vec{u_0}$. 
The initial data $\vec{u}_0$ is given by the values of $\rho_{k_0}$, 
$\cdots$, $\rho_{k_0 + M}$ in the explicit formula (\ref{p18}). 

\section{Numerical Experiments}
\setcounter{equation}{0}
The forward and inverse transforms are implemented with the 256 point FFT. For speech 
signals, Table 1 is used at sampling frequency 16 kHz. For music signals, 
Table 2 is used at sampling frequency 44.1 kHz. 
Top (bottom) panel of Figure 1 shows the 
oscillatory unpredicatibility measure $c(k)$ of a speech (music) frame. 
Top (bottom) panel of Figure 2 is the corresponding 
weighted unpredicatibility measure 
$ec(b)$ for the speech (music) frame, oscillation is slower over the coarser 
scale $b$. In Figure 3 (Figure 4), 
we compare the original and reconstructed 
FFT amplitude spectra ($k \in [20, 128]$) of a speech (music) frame. 
The difference is negligible for $k \in [0,20]$. We see that the reconstructed 
FFT spectra captured well the upper envelope of the original FFT spectra of the 
speech frame. For the music frame, much more details of the FFT spectra are recovered. 
Except for a mismatched peak and a valley over $k \in [20,40]$, the dashed and 
solid curves nearly agree. If one zooms in further, one may see differences 
over smaller scales yet the reconstructed (dashed) curve again keeps track of the 
envelope of the original spectral shape well. 
Figure 5 compares the smoother spectral solution ($\tau=2$, dashed) with 
the simple solution ($\tau=0$, solid) in case of a speech frame over $k \in [30,128]$ 
where constrained optimization (smoothing) takes place. The steady state 
is almost approached at $\tau =2$. The smoothing is 
similar for music frames. 

Figure 6 (Figure 7) compares the original and reconstructed 
speech (music) waveforms. The total relative $l^2$ error for the speech  
signal in Figure 6 is 12 \%, and is only 1.5\% 
for music signal of Figure 7. This is consistent with 
the better spectral fit of Figure 4 than that of Figure 3.  
The improvement by the optimization (\ref{p19})-(\ref{p21}) is however found to be 
minor both in terms of the relative $l^2$ error of reconstructed signals 
and perceptual difference in hearing the signals. The optimization step may 
be helpful however in other signal processing tasks to be evaluated in the future.

The original and reconstructed $(\tau=0)$ speech (music) signals 
in Figure 6 and 
Figure 7 can be heard 
at http://math.uci.edu/$\sim $jxin/sounds.html.
Inspite of the errors (loss) incurred in the reconstruction, 
there is very little perceptual difference 
between the original and the reconstructed signals, thanks to the masking effects 
present in the human ears \cite{Sch}. 
Hence we have achieved the perceptually equivalent inversion of the many-to-one 
transform. 

\begin{figure}[p]
\centerline{\includegraphics[width=350pt,height=350pt]{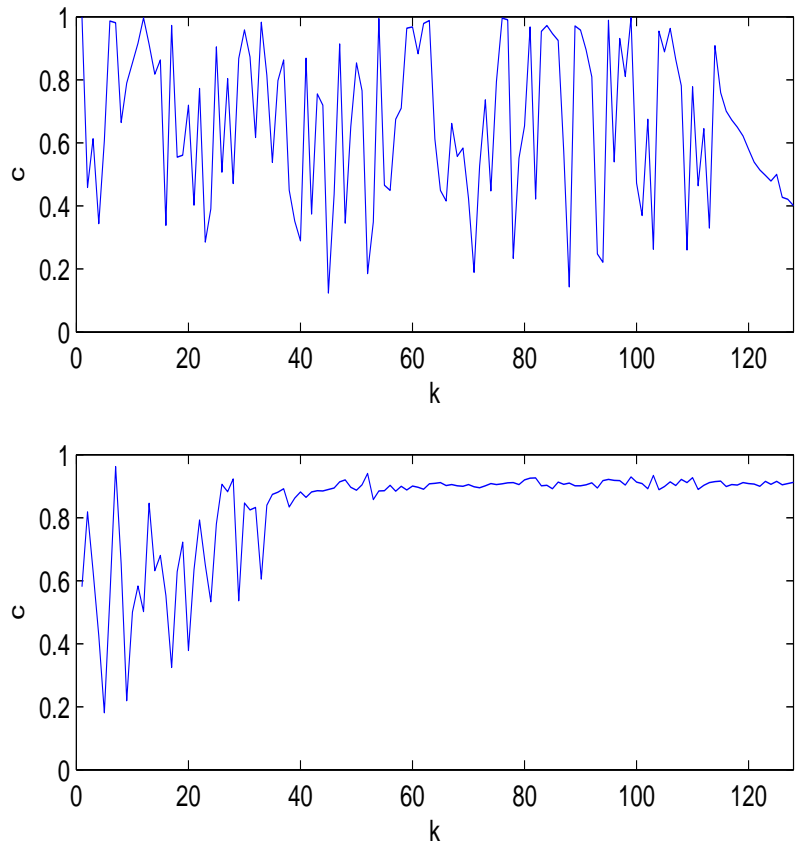}}
\vspace{.1 in}
\caption{Top panel: unpredicatibility measure $c(k)$ of a speech frame, illustrating its 
oscillatory nature in FFT wave number $k$, $k \in [0,128]$.
Bottom panel: unpredicatibility measure $c(k)$ of a music frame, $k \in [0,128]$.}
\end{figure}

\begin{figure}[p]
\centerline{\includegraphics[width=350pt,height=350pt]{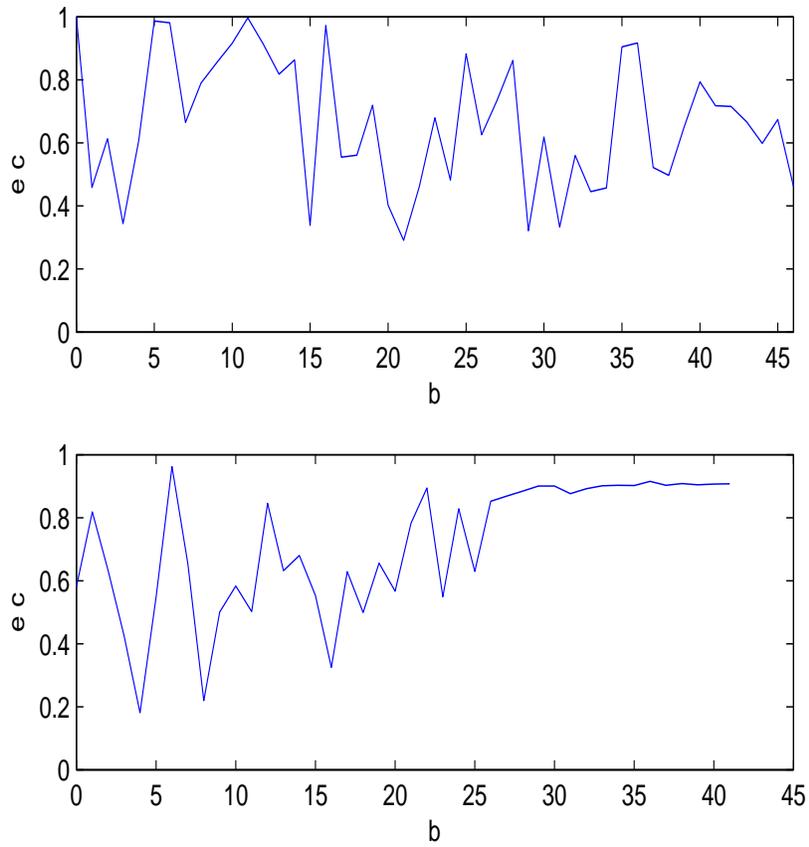}}
\vspace{.1 in}
\caption{Top panel: weighted unpredicatibility measure 
$ec(b)$ of a speech frame, $b \in [0,46]$.
Bottom panel: weighted unpredicatibility measure $ec(b)$ of a music frame, $b \in [0,41]$.}
\end{figure}

\begin{figure}[p]
\centerline{\includegraphics[width=350pt,height=350pt]{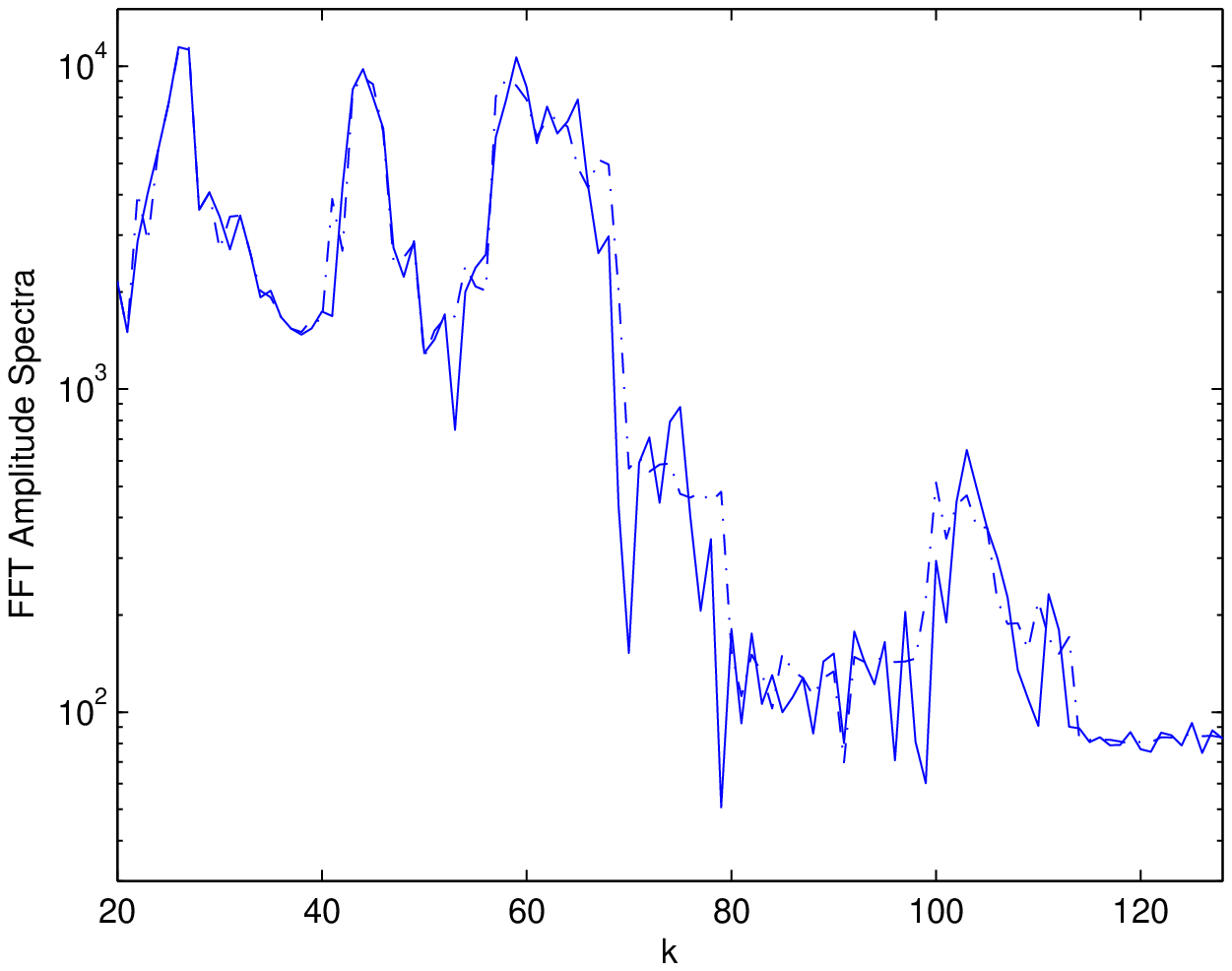}}
\vspace{.1 in}
\label{Fig5}
\caption{Original (solid) and reconstructed (dashed, $\tau=0$) FFT amplitude 
spectra of a speech frame.}
\end{figure}

\begin{figure}[p]
\centerline{\includegraphics[width=350pt,height=350pt]{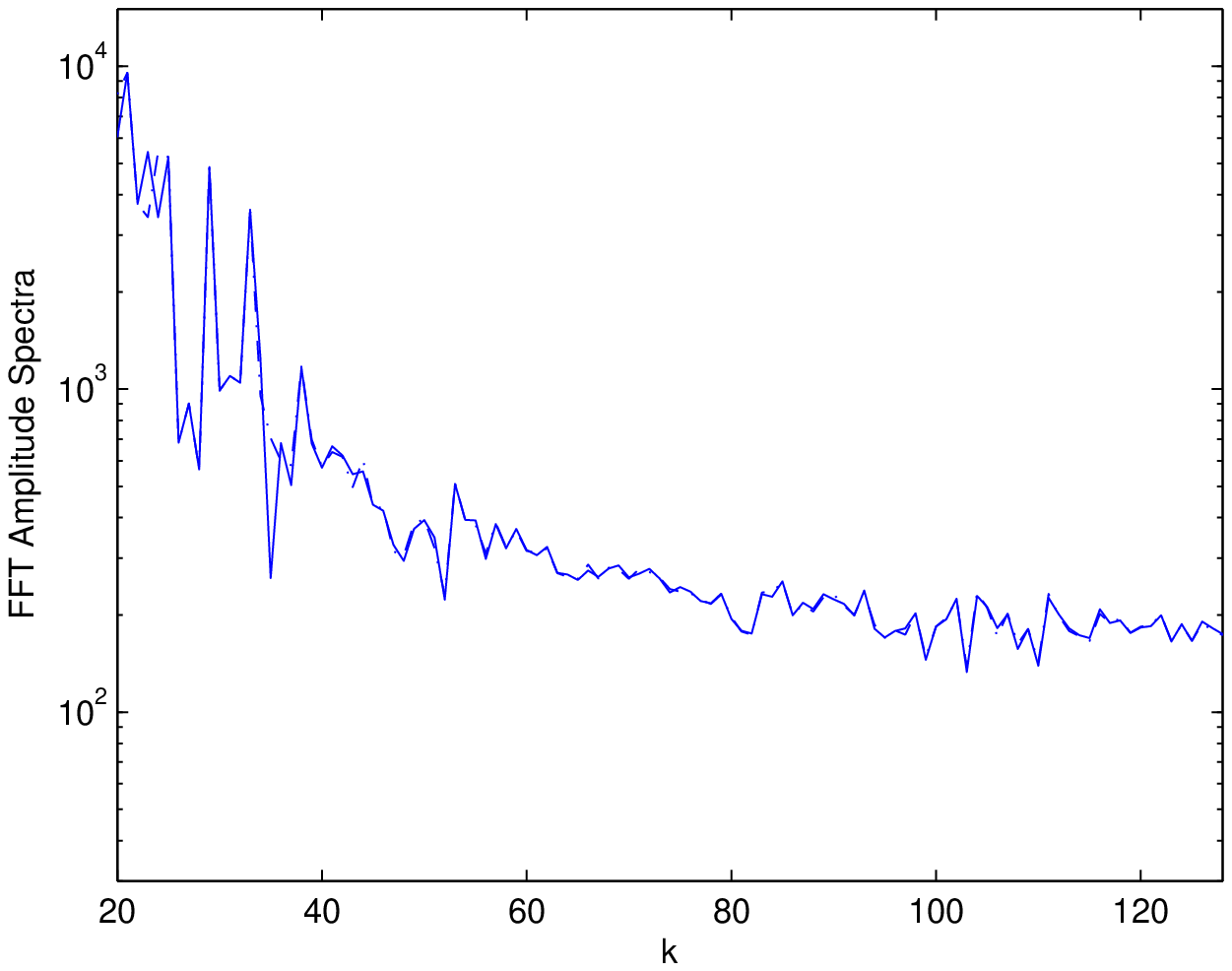}}
\vspace{.1 in}
\label{Fig6}
\caption{Original (solid) and reconstructed (dashed, $\tau=0$) FFT amplitude 
spectra of a music frame.}
\end{figure}

\begin{figure}[p]
\centerline{\includegraphics[width=350pt,height=350pt]{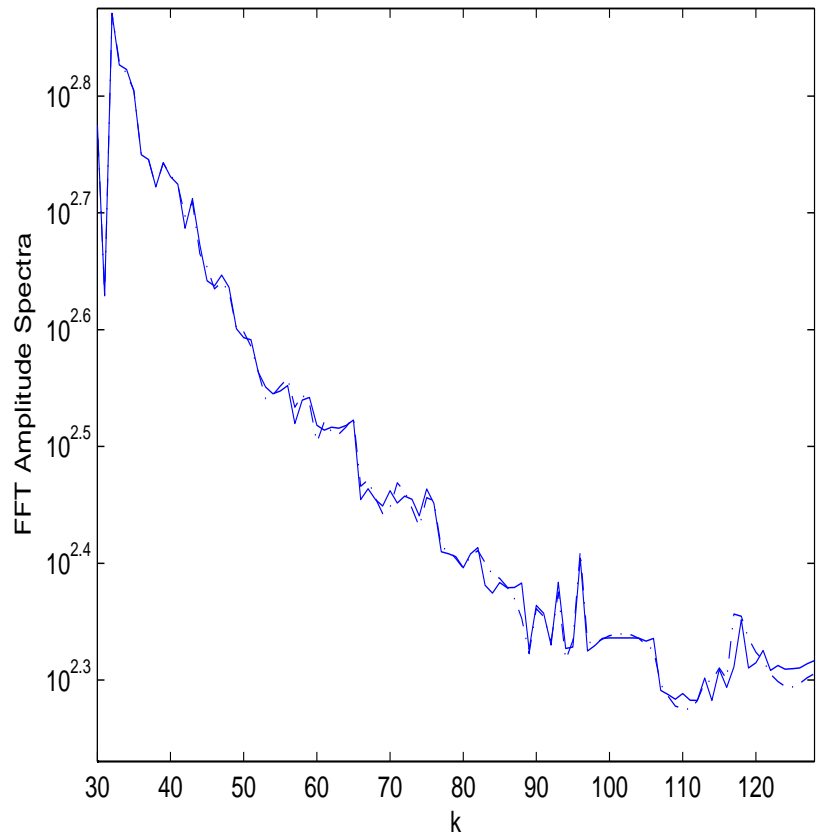}}
\vspace{.1 in}
\label{Fig7}
\caption{Comparison of reconstructed (solid, $\tau=0$) and (dashed, $\tau=2$) 
FFT amplitude spectra of a speech frame. The 
dashed curve is smoother while satisfying the same spectral constraints.}
\end{figure}

\begin{figure}[p]
\centerline{\includegraphics[width=350pt,height=350pt]{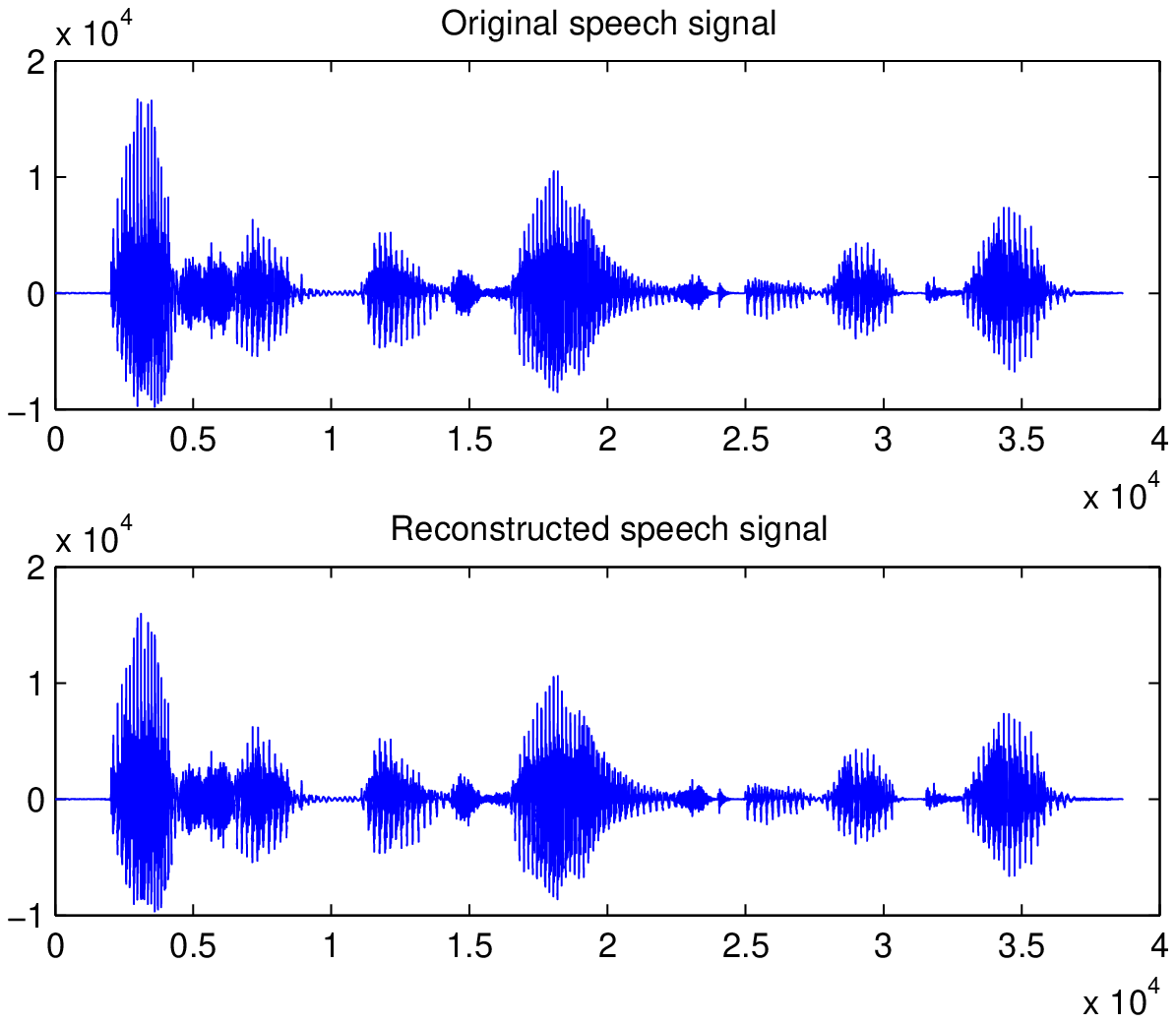}}
\vspace{.1 in}
\label{Fig8}
\caption{Comparison of the input (top) and reconstructed (bottom, $\tau=0$) speech signals in waveforms.}
\end{figure}

\begin{figure}[p]
\centerline{\includegraphics[width=350pt,height=350pt]{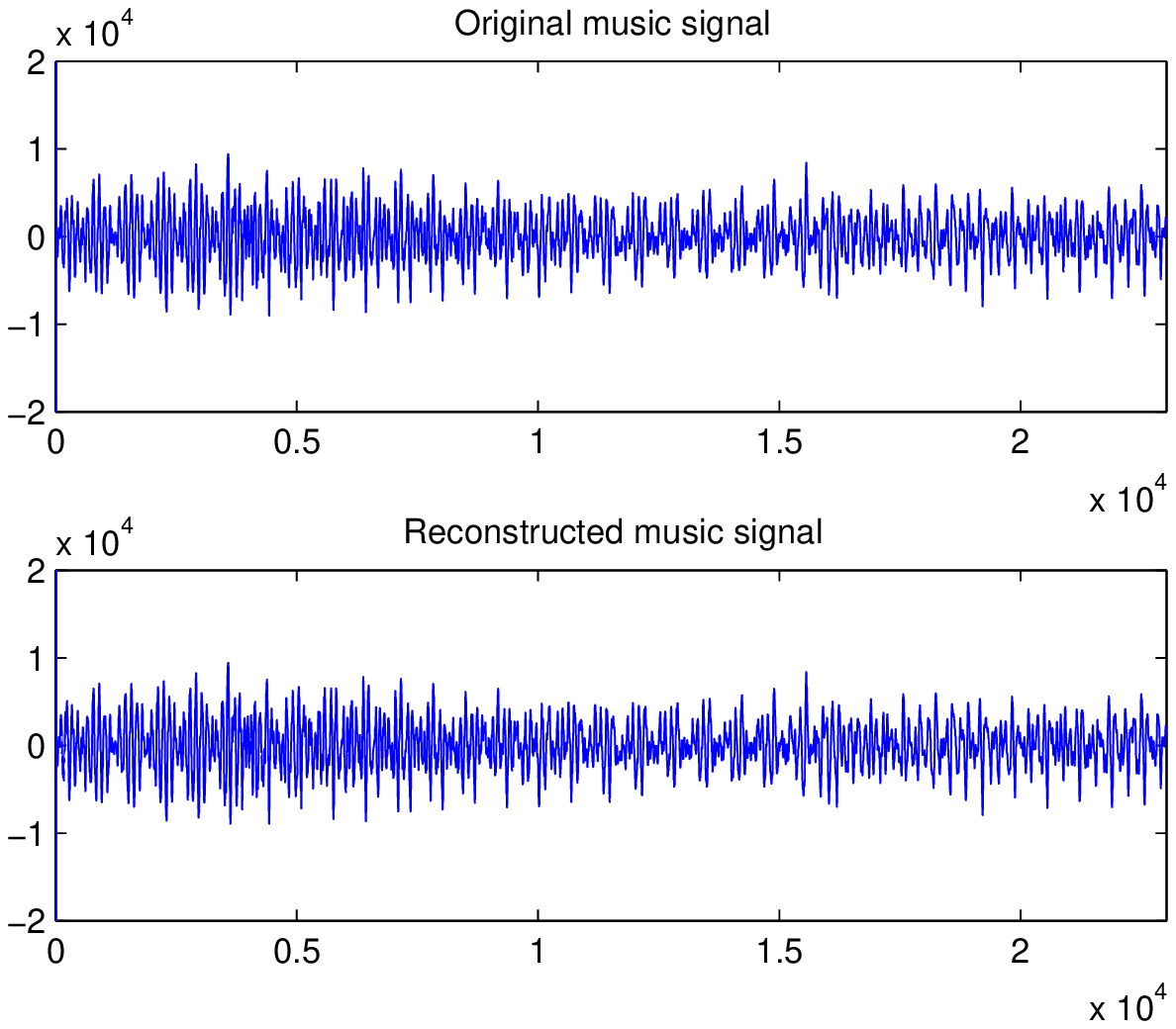}}
\vspace{.1 in}
\label{Fig9}
\caption{Comparison of the input (top) and reconstructed (bottom, $\tau=0$) music signals in waveforms.}
\end{figure}

\section{Discussion and Conclusion}
A many-to-one auditory transform is introduced so that the resulting spectrum, 
especially towards the higher frequency regime, is much less refined
than the FFT spectrum, yet just enough to resolve the band widths 
of human auditory filters (critical bands). A reconstruction of perceptually 
equivalent inverse is given so that the inverted signal 
makes little perceptual difference from 
the input signal even though there is a loss mathematically. The inversion preserves 
the band energies and band signal to noise ratios, which prove to be essential 
in capturing the perception of sounds. Both the forward and inverse transforms 
are in closed analytical form and can be carried out in real time. 
Test examples on speech and music signals  
illustrated the properties of the transform and its inversion. 
The transform is a promising new tool 
for sound compensation or enhancement that requires spectral manipulations 
over the scale of critical bands. 

A future study may concern with more accurate
inversion while conserving additional spectral information 
of the signal, such as energy variation about its mean value $e(b)/N_b$ 
inside each frequency band with $N_b \geq 3$. Another is to further develop MDAT in 
specific applications such as hearing aids and hearing implants.

\section{Acknowledgements}
This work was supported in part by NSF grant ITR-0219004 (J.X), and 
NIH grant 2R43DC005678-02A1 (Y.Q).
We thank Profs. G. Papanicolaou and H-K Zhao for helpful conversations. 

\vspace{.1 in}


\bibliographystyle{plain}

\newpage
\begin{table}
\caption{Partition and psychoacoustic parameters for the 256 point FFT at 16 kHz sampling frequency. 
The columns are (from left to right) band index, low FFT index of the band, high
FFT index of the band, number of FFT components in the band (width), the bark value of the band. 
The symmetric part of the AC components of FFT are not listed. Zero index 
refers to DC component of FFT.}
\begin{center}
\begin{tabular}{||c|c|c|c|c||}\hline
Band Index & Low FFT Index & High FFT Index & Width  &  Bark Value \\
\hline \hline
0 &  0  &  0  &  1  &  0 \\
\hline
1    &  1     &  1      & 1         &  0.63    \\
\hline
2       &  2  &  2      & 1           & 1.26    \\
\hline 
3  &     3    &  3     &   1          & 1.88     \\
\hline
4 & 4         & 4      & 1         & 2.50        \\
\hline
5   & 5        &  5    & 1          & 3.11        \\
\hline
6       &  6      & 6     &  1    &    3.70        \\
\hline
7 & 7  & 7  & 1 & 4.28 \\
\hline
8 & 8 & 8 & 1 & 4.85 \\
\hline
9  & 9  & 9  & 1 & 5.39 \\
\hline
10 & 10  & 10 & 1 & 5.92 \\
\hline 
11 & 11 & 11 & 1 & 6.43 \\
\hline
12 & 12 & 12 & 1 & 6.93 \\
\hline
13 & 13 & 13 & 1 & 7.40 \\
\hline
14 & 14 & 14 & 1 & 7.85 \\
\hline
15 & 15 & 15 & 1 & 8.29 \\
\hline
16 & 16 & 16 & 1 & 8.70 \\
\hline
17 & 17 & 17 & 1 & 9.10 \\
\hline
18 & 18 & 18 & 1 & 9.49 \\
\hline
19 & 19 & 19 & 1 & 9.85 \\
\hline
20 & 20 & 20 & 1 & 10.20 \\
\hline
21 & 21 & 22 & 2  & 10.85 \\
\hline
22 & 23 & 24 &  2 & 11.44 \\
\hline
23 & 25 & 26 & 2 & 11.99 \\
\hline
24 & 27  & 28 & 2 & 12.50 \\
\hline
25 & 29 & 30 & 2 & 12.96 \\
\hline
26 & 31 & 32 & 2 & 13.39 \\
\hline
27 & 33 & 34 & 2 & 13.78 \\
\hline
\end{tabular}\end{center}
\end{table}

\newpage
\begin{table}
\setcounter{table}{0}
\begin{center}
\caption{ Continued.}
\medskip
\begin{tabular}{||c|c|c|c|c||}\hline
Band Index & Low FFT Index & High FFT Index & Width  &  Bark Value \\
\hline \hline
28 & 35 & 36 & 2 & 14.15 \\
\hline
29 & 37 & 39 & 3 & 14.57 \\
\hline
30 & 40 & 42 & 3 & 15.03 \\
\hline
31 & 43 & 45 & 3 & 15.45 \\
\hline
32 & 46 & 48 & 3 & 15.84 \\
\hline
33 & 49 & 51 & 3 & 16.19 \\
\hline
34 & 52 & 55 & 4 & 16.57 \\
\hline
35 & 56 & 59 & 4 & 16.97 \\
\hline
36 & 60 & 63 & 4 & 17.33 \\
\hline
37 & 64 & 68 & 5 & 17.71 \\
\hline
38 & 69 & 73 & 5 & 18.09 \\
\hline
39 & 74 & 78 & 5 & 18.44 \\
\hline 
40 & 79 & 84 & 6 & 18.80 \\
\hline
41 & 85 & 90 & 6 & 19.17 \\
\hline
42 & 91 & 97 & 7 & 19.53 \\
\hline
43 & 98 & 104 & 7 & 19.89 \\
\hline
44 & 105 & 112 & 8 & 20.25 \\
\hline
45 & 113 & 120 & 8 & 20.61 \\
\hline
46 & 121 & 127 & 7 & 20.92\\ 
\hline
 \end{tabular}\end{center}
\end{table}

\newpage

\begin{table}
\caption{Partition and psychoacoustic parameters for the 256 point FFT at 44.1 kHz sampling frequency. 
The columns are (from left to right) band index, low FFT index of the band, high
FFT index of the band, number of FFT components in the band (width), 
the bark value of the band. 
The symmetric part of the AC components of FFT are not listed. Zero index 
refers to DC component of FFT.}
\begin{center}
\begin{tabular}{||c|c|c|c|c||}\hline
Band Index & Low FFT Index & High FFT Index & Width  &  Bark Value \\
\hline \hline
0 &  0  &  0  &  1  &  0 \\
\hline
1    &  1     &  1      & 1         &  1.73    \\
\hline
2       &  2  &  2      & 1           & 3.41    \\
\hline 
3  &     3    &  3     &   1          & 4.99     \\
\hline
4 & 4         & 4      & 1         & 6.45        \\
\hline
5   & 5        &  5    & 1          & 7.75        \\
\hline
6       &  6      & 6     &  1    &    8.92        \\
\hline
7 & 7  & 7  & 1 & 9.96 \\
\hline
8 & 8 & 8 & 1 & 10.87 \\
\hline
9  & 9  & 9  & 1 & 11.68 \\
\hline
10 & 10  & 10 & 1 & 12.39 \\
\hline 
11 & 11 & 11 & 1 & 13.03 \\
\hline
12 & 12 & 12 & 1 & 13.61 \\
\hline
13 & 13 & 13 & 1 & 14.12 \\
\hline
14 & 14 & 14 & 1 & 14.59 \\
\hline
15 & 15 & 15 & 1 & 15.01 \\
\hline
16 & 16 & 16 & 1 & 15.40 \\
\hline
17 & 17 & 17 & 1 & 15.76 \\
\hline
18 & 18 & 19 & 2 & 16.39 \\
\hline
19 & 20 & 21 & 2 & 16.95 \\
\hline
20 & 22 & 23 & 2 & 17.45 \\
\hline
21 & 24 & 25 & 2  & 17.89 \\
\hline
22 & 26 & 27 &  2 & 18.30 \\
\hline
23 & 28 & 29 & 2 & 18.67 \\
\hline
24 & 30  & 31 & 2 & 19.02 \\
\hline
25 & 32 & 34 & 3 & 19.41 \\
\hline
26 & 35 & 37 & 3 & 19.85 \\
\hline
27 & 38 & 40 & 3 & 20.25 \\
\hline
\end{tabular}\end{center}
\end{table}

\newpage
\begin{table}
\setcounter{table}{1}
\begin{center}
\caption{ Continued.}
\medskip
\begin{tabular}{||c|c|c|c|c||}\hline
Band Index & Low FFT Index & High FFT Index & Width  &  Bark Value \\
\hline \hline
28 & 41 & 43 & 3 & 20.62 \\
\hline
29 & 44 & 47 & 4 & 21.01 \\
\hline
30 & 48 & 51 & 4 & 21.43 \\
\hline
31 & 52 & 55 & 4 & 21.81 \\
\hline
32 & 56 & 59 & 4 & 22.15 \\
\hline
33 & 60 & 64 & 5 & 22.51 \\
\hline
34 & 65 & 69 & 5 & 22.87 \\
\hline
35 & 70 & 75 & 6 & 23.23 \\
\hline
36 & 76 & 81 & 6 & 23.59 \\
\hline
37 & 82 & 88 & 7 & 23.93 \\
\hline
38 & 89 & 96 & 8 & 24.00 \\
\hline
39 & 97 & 105 & 9 & 24.00 \\
\hline 
40 & 106 & 115 & 10 & 24.00 \\
\hline
41 & 116 & 127 & 12 & 24.00 \\
\hline
\end{tabular}\end{center}
\end{table}

\newpage









\end{document}